\documentclass[12pt,oneside,english]{amsart}
\usepackage[latin1]{inputenc}
\usepackage{amssymb}

\makeatletter
\newtheorem{theorem}{Theorem}

\newtheorem{example}{Example}

\newtheorem{remark}{Remark}
\newtheorem{conjecture}{Conjecture}
\newtheorem{question}{Question}
\newtheorem{corollary}{Corollary}

\usepackage{babel}

\makeatother
\begin{document}
\baselineskip=17pt

\title[On a unique additive representations of integers ]{On unique additive representations of positive integers and some close problems}

\author{Vladimir Shevelev}
\address{Departments of Mathematics \\Ben-Gurion University of the
 Negev\\Beer-Sheva 84105, Israel. e-mail:shevelev@bgu.ac.il}


\begin{abstract}
Let, for $r\geq 2,$ $\{m_r(n)\}_{n\geq0}$ be Moser sequence such that every nonnegative integer is the unique sum of the form $m_r(k)+rm_r(l).$ In this article we give an explicit decomposition formulas of such form and a surprisingly simple recursion relations for Moser's numbers. We also study interesting properties of the sequence $\{rm_r(n-1)+1\}_{n\geq1}$ and its connection with some important problems. In particular, in case of $r=2,$ this sequence is surprisingly connected with the numbers solving the combinatorial Josephus-Gro\"{e}r  problem. We pose also some open questions.

\end{abstract}

\maketitle

\section{Introduction}

     As well known, I. M. Vinogradov in 1937 proved that every sufficiently large odd $n$ is represented as a form $n=p_1+p_2+p_3,$ where $p_i$ are primes, i. e. he solved with not more than finite set of exceptions the so-called ternary Goldbach conjecture. This problem is essentially weaker than so-called binary Goldbach conjecture, which states that every even number is represented as a form  $n=p_1+p_2.$ Another natural strengthening of the ternary
     problem is the following very interesting conjecture.
     \begin{conjecture}\label{1}(Lei Zhou (2005), see [7], comment to sequence A103151)
Every odd number more or equal than 9 is represented as a form $p_1+2p_2,$ where p and q are odd primes.
\end{conjecture}
As the Goldbach binary representation, the representation of Conjecture 1, generally speaking, is not unique. Let $B$ be the set of analytical functions $f(x)\enskip |x|<1$, with (0,1)-Taylor coefficients. In 1962, L. Moser ([6], see also [4]) showed that the functional equation
\begin{equation}\label{1}
f(x)f(x^r)=\frac {1} {1-x}, \enskip |x|<1.
\end{equation}
has the unique solution in $B$
\begin{equation}\label{2}
 f(x)=\sum^{\infty }_{n=0}x^{m_r(n)}, \enskip m_r(n)\geq 1, \enskip |x|<1.
\end{equation}
This means that for every $r\geq2$ there exists the unique sequence $\{m_r(n)\}_{n\geq0}$ such that every $n\geq0$ could be represented by a unique way as a sum:\newpage
\begin{equation}\label{3}
n=m_r(k)+rm_r(l)
\end{equation}
 In the case of $r=2,$ the sequence also was studied by N. G. de Bruijn [3] in connection with the representations of positive integers by differences of the form $am_2(k)-bm_2(l)$ with positive coefficients $a, b.$  Therefore, in this special case, it is known as Moser-de Bruijn sequence (see also [1, example 13] and [7, sequence A000695]). Its the first terms are:
\begin{equation}\label{4}
0, 1, 4, 5, 16, 17, 20, 21, 64, 65, 68, 69, 80, 81, 84, 85, 256, 257,...
\end{equation}
 Note that, using Moser's method, one can prove that, for a given odd $c\geq3,$  every odd number $N\geq3c$ is  unique sum of the form $a^{(c)}(k)+2a^{(c)}(l)$, where
 \begin{equation}\label{5}
 a^{(c)}(n)=2m_2(n)+c,
 \end{equation}
 (in particular, in case of $c=3$ see sequence A147568 in [7]); moreover, this sequence is unique sequence with such property. Now we prove the following.
\begin{theorem}\label{1}
For every odd $c$, the sequence (5) contains infinitely many composite numbers.
\end{theorem}
\bfseries Proof.\enskip\mdseries As is known [1, example 13], the terms of sequence (4) are sums of distinct powers of 4. Therefore, e.g., pairs of the form  $(2\cdot4^k, \enskip 2\cdot4^k+2) $ and $(2\cdot4^k+2\cdot4^l,\enskip 2\cdot4^k+2\cdot4^l+2)$  are terms of $\{2m_2(n)\}.$ Since, furthermore, $ 4^k\equiv1\pmod3,$ then $2\cdot4^k+1\equiv3\pmod6,$  while $2\cdot4^k+2\cdot4^l+1\equiv5\mod6,$ i.e. we have infinitely many such numbers with different residues modulo 6.
Therefore, for arbitrary given $c$, the sequence (5) contains  infinitely many pairs with difference 2,
such that middle numbers have different residues modulo 6. This means that these pairs contain infinitely many composite numbers.$\blacksquare$\newline
Hence, in connection with Conjecture 1, we obtain the following statement.
\begin{corollary}
 There is no a subsequence $T$ of primes, such that every odd number $ m\geq C$ is expressible as unique sum of the form $m=p+2q,\enskip p$ and $q $ are in $T.$
\end{corollary}
\begin{conjecture}\label{2}
 For arbitrary positive numbers $C_1$ and $C_2$, there is no a subsequence $T$ of primes, such that every odd number $ m\geq C_1$ is expressible by the number $\nu,\enskip 1\leq\nu\leq C_2,$ of different ways as a sum of the form $m=p+2q,\enskip p$ and $q $ are in $T.$
 \end{conjecture}
 \begin{conjecture}\label{3}
 For every odd $c$, the sequence $\{a^{(c)}(n)\}$ (5) contains infinitely many primes.
 \end{conjecture}\newpage
Furthermore, using the method of generating functions, it is easy to show that the function
 \begin{equation}\label{6}
 g(x)=\sum^{\infty }_{n=1}x^{s_n^{(r)}}, \enskip{s_n^{(r)}}\geq 1, \enskip |x|<1,
 \end{equation}
 where
  \begin{equation}\label{7}
s_n^{(r)}=rm_r(n-1)+1, \enskip n\geq1,
  \end{equation}
 is the unique solution in $B$ of the functional equation
 \begin{equation}\label{8}
 f(x)f(x^r)=\frac {x^{r+1}} {1-x^r}, \enskip |x|<1.
 \end{equation}
 This means that for every $r\geq2$  there exists the unique sequence $\{s_n^{(r)}\}$ such that every $n\equiv1 \pmod {r}, \enskip n\geq r+1,$ could be represented by a unique way as a sum:
\begin{equation}\label{9}
n=s_k^{(r)}+rs_l^{(r)}.
\end{equation}
In this article we give an explicit decomposition formulas of the form (3) and (9) and an unexpectedly simple recursion relations for Moser's numbers $m_r(n).$  We also study very interesting properties of the sequence $\{s_n^{(r)}\}$ and its connection with some important problems. In particular, in case of $r=2,$ the sequence  $\{s_n^{(2)}\}$  is surprisingly connected with the numbers solving the combinatorial Josephus-Gro\"{e}r  problem (see [5] and sequences A088442, A090569 in [7]). We also pose several questions.
\section{Explicit decomposition formulas}
\begin{theorem}\label{2}
 Let $n\geq 0$ and
$$n=\sum_{i\geq0}\nu_ir^{i},\enskip 0\leq\nu_i\leq r-1,$$
be the expansion of  $n$ in base $r\geq2,$ then\newline
1)$$ m_r(n)=\sum_{i\geq0}\nu_ir^{2i};$$
2)the diophantine equation
$$m_r(k)+rm_r(l)=n$$
has the unique solution
$$k=\sum_{i\geq 0}\nu_{2i}r^{i}, \enskip l=\sum_{i\geq 0}\nu_{2i+1}r^{i}.$$
\end{theorem}
\begin{example}\label{1}
Let $n=27, r=2.$ Then $\nu_0=1, \nu_1=1, \nu_2=0, \nu_3=1, \nu_4=1.$\newpage Therefore, we have
$$ m_2(27)=4^4+4^3+4+1=325;$$
$$k=\nu_0+\nu_2\cdot2+\nu_4\cdot2^2=5$$
and
$$l=\nu_1+\nu_3\cdot2=3,$$
such that
$$m_2(5)=17, \enskip m_2(3)=5\enskip and \enskip 27=17+2\cdot5.$$
\end{example}
Proof of Theorem 2 is similar to one of the following theorem.
\begin{theorem}\label{3}
1) Let $n\geq 1.$ If
$$n-1=\sum_{i\geq0}\nu_ir^{i},\enskip 0\leq\nu_i\leq r-1,$$
is the expansion of  $n-1$ in base $r,$ then
$$s^{(r)}_n=1+\sum_{i\geq0}\nu_ir^{2i+1};$$
2) If, for given $r\geq2,\enskip N\equiv 1 \pmod {r},\enskip N\geq r+1,$ the number $N-r$ has the expansion in base $r$
$$N-r=\sum_{j=0,2,4,...}a_jr^{j}+\sum_{j=1,3,5,...}b_jr^{j},$$
then the diophantine equation
$$s^{(r)}_k+rs^{(r)}_l=N$$
has the unique solution
$$k=1+\sum_{i\geq 0}b_{2i+1}r^{i},$$
$$l=1+\sum_{i\geq 0}a_{2i+2}r^{i}.$$
\end{theorem}
\bfseries Proof.\enskip\mdseries By the expansion of $N-r$, where according to the condition $a_0=1,$ we have

$$N=(r-1+\sum_{j=0,2,4,...}a_jr^{j})+(1+\sum_{j=1,3,5,...}b_jr^{j})=$$
$$r(1+\sum_{j=2,4,...}a_jr^{j-1})+(1+\sum_{j=1,3,5,...}b_jr^{j})=$$
\begin{equation}\label{10}
r(1+\sum_{j=1,3,5,...}a_{j+1}r^{j})+(1+\sum_{j=1,3,5,...}b_jr^{j}).
\end{equation}
Let us consider a set $E^{(r)}$ of positive integers $u$ for which the expansion in base $r$ has the form:\newpage
  \begin{equation}\label{11}
u=1+\sum_{j=1,3,5,...}d_jr^{j}.
\end{equation}
Denote $\{u_n\}$ the sequence which is obtained after the lexicographical ordering of set $E^{(r)}.$  Then  $(u_n)$ is an increasing sequence and  in view of one-to-one lexicographical correspondence $n-1\leftrightarrow u_n$ (such that $n-1 \leftrightarrow u_1=1$) we conclude that
if
\begin{equation}\label{12}
n-1=\sum_{i\geq 0}\nu_ir^{i},
\end{equation}
 then
 \begin{equation}\label{13}
 u_n=1+\sum_{i\geq 0}\nu_{i}r^{2i+1},
\end{equation}
Thus, (10) shows that for every positive integer $N\equiv1 \pmod r$ there exist positive integers $k,l$ such that
$$N=u_k+ru_l.$$
  To demonstrate the uniqueness of this representations, put
$$N=(1+\sum_{j=1,3,5,...}g_jr^{j})+r(1+\sum_{j=1,3,5,...}h_{j}r^{j}).$$
Then we have
$$N=1+\sum_{j=1,3,5,...}g_jr^{j} + r-1 + (1+\sum_{j=1,3,5,...}h_{j}r^{j+1})=$$
$$r+\sum_{j=1,3,5,...}g_jr^{j}+\sum_{j=0,2,4,...}h_{j-1}r^{j},$$
where $h_{-1}=1.$ Hence,
$$N-r=\sum_{j=1,3,5,...}g_jr^{j}+\sum_{j=0,2,4,...}h_{j-1}r^{j}$$
and the comparison with the condition of 2) shows that $g_j=b_j, h_{j}=a_{j+1},$ which corresponds to (10). In view of the uniqueness of such sequence, we conclude that $s_n=u_n,$ and from (13) the statement 1) of the theorem follows.
Furthermore, in view of decomposition (10), we have
$$s^{(r)}_k=1+\sum_{i\geq0}b_{2i+1}r^{2i+1},$$
$$s^{(r)}_l=1+\sum_{i\geq0}a_{2i+2}r^{2i+1}.$$\newpage
Thus, by the one-to-one lexicographical correspondence $n-1\leftrightarrow u_n$ defined by (12)-(13), we obtain the required expressions for $k$ and $l.\blacksquare$
\begin{remark}
Let $A^{(r)}(x)$ be the counting function of the terms $s^{(r)}_n$ not exceeding $x$. It could be proved that for every $r\geq2,\enskip A^{(r)}(x)=O(\sqrt{x})$ and $A^{(r)}(x)=\Omega(\sqrt{x}).$
\end{remark}
\section{Fast calculation of consecutive values of $ m(n)=m_r(n) \enskip and\enskip s_n=s^{(r)}_n  $}
Now we prove a simple recursion relation for sequences $m(n)=m_r(n)$ and $s_n=s^{(r)}_n,$ taking into
account that, by (1),(2) and (7),for every $r\geq2,\enskip m_r(0)=0$ and $s_1^{(r)}=1.$
\begin{theorem}\label{4}
For $n\geq1,$
\begin{equation}\label{14}
 m(n)=m(n-1)+\frac {1} {r+1}(r^{2\tau+1}+1),
 \end{equation}
 where $\tau\geq0$ is such that $n\equiv r^\tau\pmod {r^{\tau+1}}$;\newline
for $n\geq2,$
\begin{equation}\label{15}
s_n=s_{n-1}+\frac {1} {r+1}(r^{2(t+1)}+r),
\end{equation}
where $t\geq0$ is such that $n-1\equiv r^t\pmod {r^{t+1}}.$
\end{theorem}
\bfseries Proof.\enskip\mdseries In view of (7), it is sufficient to prove only one of these formulas.
Below we prove the second of them. \newline Let
$$n=\sum_{i\geq 0}b_ir^{i}, \enskip 0\leq b_i\leq r-1$$
be the expansion of $n$ in base $r$ and $j\geq0$ be the least index for which $b_j>0,$ such that
 $$n=b_jr^j+\sum_{i\geq j+1}b_ir^{i}$$
 and
 $$n-1=(b_j-1)r^j+(r-1)(1+r+...+r^{j-1})+\sum_{i\geq j+1}b_ir^{i}.$$
 Then, by (12)-(13), we have
\begin{equation}\label{16}
s_n=1+(b_j-1)r^{2j+1}+(r-1)(r+r^3+...+r^{2j-1})+\sum_{i\geq j+1}b_ir^{2i+1}.
\end{equation}
1)If here $j\geq1,$ then
$$n-2=(b_j-1)r^j+(r-2)+(r-1)(r+...+r^{j-1})+\sum_{i\geq j+1}b_ir^{i}$$\newpage
and  by (12)-(13), we find
$$s_{n-1}=1+(b_j-1)r^{2j+1}+(r-2)r+(r-1)(r^3+...+r^{2j-1})+\sum_{i\geq j+1}b_ir^{2i+1}.$$
Thus, in view of (16), in case of $n\equiv0\pmod r$ we have
$$s_n-s_{n-1}=(r-1)r-(r_2)r=r,$$
which corresponds to the theorem in case of $t=0.$\newline
2a) Let $j=0$ \enskip and \enskip $b_0\geq2.$
Then
$$n-1=b_0-1+\sum_{i\geq1}b_ir^{i},\enskip n-2=b_0-2+\sum_{i\geq1}b_ir^{i}$$
and,consequently,
$$s_n=1+(b_0-1)r+\sum_{i\geq1}b_ir^{2i+1},\enskip s_{n-1}=1+(b_0-2)r+\sum_{i\geq1}b_ir^{2i+1}.$$
Thus, again in this case we have
$$s_n-s_{n-1}=(r-1)r-(r-2)r=r.$$
This completes the case $t=0$ of the theorem.\newline
2b) At last, let $j=0,$ \enskip $b_0=1,$ \enskip i. e.
$$n=1+ \sum_{i\geq1}b_ir^{i}.$$
Let, furthermore, $t\geq1$ be the least index for which $b_t>0,$ such that
 $$n-1=b_tr^t+\sum_{i\geq t+1}b_ir^{i}.$$
 Then
 $$n-2=(b_t-1)r^t+(r-1)(1+r+...+r^{t-1})+\sum_{i\geq t+1}b_ir^{i}.$$
 and we find
 $$s_{n-1}=1+(b_t-1)r^{2t+1}+(r-1)(r+r^3+...+r^{2t-1})+\sum_{i\geq t+1}b_ir^{2i+1},$$
 while
 $$s_n=1+b_tr^{2t+1}+\sum_{i\geq t+1}b_ir^{2i+1},$$
 and after simple transformation we obtain the theorem. $\blacksquare$\newline
 In particular, in case of $r=2,$ for calculation of the consecutive values of sequence $\{s^{(2)}_n\}$  we have a recursion:\newpage
\begin{equation}\label{17}
s_n=s_{n-1}+\frac {1} {3}(4^{t+1}+2), \enskip n\geq2, \enskip s_1=1,
\end{equation}
where $t\geq0$ is such that $n-1\equiv 2^t\pmod {2^{t+1}}.$
\section{Case of $r=2$: connection with Josephus-Gro\"{e}r problem}
In [5] Chris Gro\"{e}r posed and investigated a modification of well known Josephus combinatorial problem.
It is interesting that the case of $r=2$ of Theorem 1 is connected with the following special case of this problem (sf [5], $q=2$):"Arrange $N$ persons {1,2,...,N} consecutively on a line rather than around in a circle. Beginning at the left end of the line, we remove every second person until we reach the end of the line. At this point we immediately reverse directions, taking care not to "double count" the person at the end of the line, and continue to eliminate every second person, but now moving right to left. We continue removing people in this back and forth manner until there remains a lone survivor $W(N)".$
Let $\sum^{k}_{j=0}b_j2^{j},$ where $b_j=0\enskip or\enskip 1, $ for $j=1,...,k-1$ and $b_0=b_k=1,$ be the binary expansion of $N$ or $N-1,$ whichever  is odd, i. e.
\begin{equation}\label{18}
\sum^{k}_{j=0}b_j2^{j}=\begin{cases} N, \enskip if \;\;N\enskip is\enskip odd,\\
N-1, \enskip if \;\;N\enskip is \enskip even.\end{cases}
\end{equation}
 Then
 \begin{equation}\label{19}
 W(N)=1+\sum^{k}_{j=1, j\enskip odd}b_j2^{j}.
 \end{equation}
Consequently, it follows from (12) that for every $N$, $W(N)\in E^{(2)}.$ Furthermore, he proved that
\begin{equation}\label{20}
W(2m)=W(2m-1)=2m+1-2W(m).
 \end{equation}
 Putting here $2m+1=N$, we see that
\begin{equation}\label{21}
N=W(N-2)+2W((N-1)/2).
\end{equation}
Since the numbers $W(N-2)$ and $W((N-1)/2)$ belong to $E^{(2)},$ i.e. are $s^{(2)}_n$-numbers, then in view of proved
in the above the uniqueness of such decomposition, we conclude that in Theorem 1 for $r=2$ we have
\begin{equation}\label{22}
 s^{(2)}_k=W(N-2),\enskip s^{(2)}_l=W((N-1)/2),
\end{equation}
i.e.,in connection with the Josephus-Gro\"{e}r problem, have the concrete combinatorial sense.\newpage
\begin{question}\label{1}
To find, if it is possible, a combinatorial sense of numbers $s^{(r)}_k$ and $s^{(r)}_l$ for $r\geq 3.$
\end{question}
\section{Some identities and iterations}
We start with the following surprising identity.
\begin{equation}\label{23}
W(s_n-2)=s_{n-1}, \enskip n\geq2.
\end{equation}
\bfseries Proof.\enskip\mdseries Indeed, let
\begin{equation}\label{24}
n-1=\sum_{i\geq 0}\nu_i2^{i}
\end{equation}
be the binary expansion of $n-1.$ Then, by (12)-(13) for $n=2,$ we have
\begin{equation}\label{25}
s_n=1+\sum_{i\geq 0}\nu_{i}2^{2i+1}.
\end{equation}
Let, furthermore, $j\geq0$ is the minimal index for which $\nu_j=1.$ Then
$$s_n-2=2^{2j+1}-1+\sum_{i\geq j+1}\nu_{i}2^{2i+1}=1+2+2^2+2^3+...+2^{2j}+\sum_{i\geq j+1}\nu_{i}2^{2i+1}.$$
Thus, by (18)-(19),
\begin{equation}\label{26}
W(s_n-2)=1+(2+2^3+...+2^{2j-1})+\sum_{i\geq j+1}\nu_{i}2^{2i+1}.
\end{equation}
On the other hand,
\begin{equation}\label{27}
n-2=2^{j}-1+\sum_{i\geq j+1}\nu_{i}2^{i}=1+2+2^2+...+2^{j-1}+\sum_{i\geq j+1}\nu_{i}2^{i}.
\end{equation}
and, by (24)-(25),
\begin{equation}\label{28}
s_{n-1}=1+(2+2^3+...+2^{2j-1}+\sum_{i\geq j+1}\nu_{i}2^{2i+1}.
\end{equation}
Now (23) follows from (26) and (28). $\blacksquare$\newline
Putting in (21) $N=s_n$ and using (23), we find
$$s_n=s_{n-1}+2W(\frac {s_n-1} {2}), \enskip n\geq2.$$
\begin{corollary}\label{2}
\begin{equation}\label{29}
\sum^{N}_{n=2}W(\frac {s_n-1} {2})=\frac {s_N-1} {2}.
\end{equation}
\end{corollary}
Besides, by the comparison with (18), we obtain the following identity:\newpage
\begin{corollary}\label{3}
\begin{equation}\label{30}
W(\frac {s_n-1} {2})=\frac {1} {6}(4^{t+1}+2), \enskip n\geq2
\end{equation}
where $t\geq0$ is such that $n-1\equiv 2^t\pmod {2^{t+1}}.$
\end{corollary}
\begin{corollary}\label{4}
Put
\begin{equation}\label{31}
V(n)=W(n-2), \enskip n\geq3.
\end{equation}
Then, beginning with every positive integer $N\geq3,$ we obtain that one of the iterations
\begin{equation}\label{32}
V(N), V^{(2)}(N)=V(V(N)), V^{(3)}=V(V^{(2)}(N)),...
\end{equation}
equals to 1.
\end{corollary}
\bfseries Proof.\enskip\mdseries Indeed, according to (22), for some $k$ we have
$$V(N)=W(N-2)=s_k.$$
Thus, by (23), the sequence (32) is
$$s_k, s_{k-1},..., s_1=1.\blacksquare$$

\section{An analog of Euler $\varphi$-function}
In the correspondence of the case of $r=2$ of Theorem 3, two odd numbers $M>1$ and $N>1$ we call collinear $(M\parallel N)$, if in decompositions $$M=s^{(2)}_{k(M)}+s^{(2)}_{l(M}),\enskip N=s^{(2)}_{k(N)}+s^{(2)}_{l(N)}$$
either $k(M)=k(N)$ or $l(M)=l(N).$ We agree that 1 is not collinear to none of $N>1.$  On the set of all positive odd numbers, consider the following function which in some sense is similar to Euler $\varphi$-function.
\begin{equation}\label{33}
\psi(1)=0,\enskip \psi(n)=\sum_{i<n, \enskip i\nparallel n}1.
\end{equation}
We prove the following property of $\psi(n).$
\begin{theorem}\label{5}
For every $m\geq1$,
$$\psi(4m-1)=\psi(4m+1).$$
\end{theorem}
\bfseries Proof.\enskip\mdseries Note that, in view of (22),
$$s^{(2)}_{l(4m-1)}=W(2m-1),\enskip s^{(2)}_{l(4m+1)}=W(2m).$$
Therefore, according to (20),
$$s^{(2)}_{l(4m-1)}=s^{(2)}_{l(4m+1)}$$
and, since sequence $\{s^{(2)}_{n}\}$ is increasing, then we conclude that\newpage
$$l(4m+1)=l(4m-1).$$
Furthermore, since
$$4m-1=s^{(2)}_{k(4m-1)}+s^{(2)}_{l(4m-1)},\enskip 4m+1=s^{(2)}_{k(4m+1)}+s^{(2)}_{l(4m+1)},$$
then
$$s^{(2)}_{k(4m+1)}=s^{(2)}_{k(4m-1)}+2.$$
Since all terms of sequence $\{s^{(2)}_{n}\}$ are odd, this could be only when
$$k(4m+1)= k(4m-1)+1.$$
Therefore, for odd integers $j<4m-1$ the number of repetitions of values $k(4m-1)$ equals to the number of repetitions of values $k(4m+1).$ From this, evidently, the theorem follows.$\blacksquare$\newline
The first terms of sequence $\{\psi(4m-1)\}$ are:
$$1,2,3,4,7,8,9,10,13,14,15,16,19,20,21,22,29,30,... $$
\begin{question}\label{2}
To find a formula for fast evaluation and asymptotics of $\psi(n).$
\end{question}
\section{Arithmetic progressions}
Only for simplicity, consider the case $r=2$. Let $a\geq1, b\geq1$ be integers. Put
\begin{equation}\label{34}
s_n(a,b)=a+b(s^{(2)}_n-1)
\end{equation}
\begin{theorem}\label{6}
If $s_n(a,b)>0,$ then for every term of the arithmetic progression
$$3a, 3a+2b, 3a+4b,...$$
there exists a unique representation of the form
\begin{equation}\label{35}
3a+2bn=s_k(a,b)+2s_l(a,b), \enskip n=0,1,2,...
\end{equation}
\end{theorem}
\bfseries Proof.\enskip\mdseries From (6), (8) in the case of $r=2$ it follows that if
$$ g(x)=\sum_{n\geq1}x^{s^{(2)}_{n}},\enskip |x|<1,$$
then
$$g(x)g(x^2)=\frac {x^3} {1-x^2}.$$
Put
$$g_{a, b}(x)=\sum_{n\geq1}x^{s_n(a,b)}=$$
$$x^{a-b}\sum_{n\geq1}x^{bs_n^{(2)}}=x^{a-b}g(x^b).$$
Then\newpage
$$g_{a, b}(x)g_{a, b}(x^2)=x^{3(a-b)}g(x^b)g(x^{2b})=\frac {x^{3a}} {1-x^{2b}}$$
and the theorem follows.$\blacksquare$

E.g., the first terms of sequence $\{s_{n}(1,2)\}$ are:
$$1,5, 17, 21, 65, 81, 85, 257,...$$
(sf sequence $A145818$ in $[7]$). Hence, every integer N of the form $4m-1, m\geq1$ has a unique representation as
$$N=s_k(1,2)+2s_l(1,2).$$

In addition, note that every odd integer
$$m\equiv3\pmod {2^{2r-1}}$$
is unique sum of the form
$$s_{2^{r-1}(k-1)+1}(1,1)+2s_{2^{r-1}(l-1)+1}(1,1),\enskip r=1,2,...,$$
while every odd integer
$$m\equiv3\pmod {2^{2r}}$$
is unique sum of the form
$$s_{2^{r-1}(k-1)+1}(1,2)+2s_{2^{r-1}(l-1)+1}(1,2),\enskip r=1,2,... .$$
\section{Binary additive representations of even numbers}
Consider the union $\{t_{n}\}$ of sequences $\{s_{n}(1,1)\}$ and $\{s_{n}(1,2)\}$ with double repetition of 1 in nondecreasing order:
$$1, 1, 3, 5, 9, 11, 17, 21, 33, 35, 41, 43, 65, 69, 81, 85, 129, 131, 137, 139, 161, 163, ... $$
(sf sequence $A145819$ in $[7]$). The following result is proved straightforward.
\begin{theorem}\label{7}
 For every even integer $N$ there exists a representation of the form
\begin{equation}\label{36}
N=t_k+t_l.
\end{equation}
\end{theorem}
Note that, equation (1) in case of $r=1,$ evidently, has not a solution in set B. Hence, it is impossible to attain the uniqueness of the binary additive representations of all even numbers. Therefore, generally speaking, we have not unique decomposition of the form (36).\newpage
\begin{conjecture}\label{4}
Let $T(x)$ be the counting function for $t_n\leq x.$ Then the sequence $\{t_{n}\}$ is minimal in the following sense: if any sequence has the counting function $ A(x)\leq T(x)$ for all $x\geq1$ and $A(x)<T(x)$ for $x\geq x_0,$ then there exists an even integer $N$ which is not expressible as a sum of two terms of such sequence.
\end{conjecture}
In the connection of Conjecture 4, it is interesting the sequence $\{v_{n}\}$ of the least even positive integers which is not expressible as a sum $t_k+t_l,$ if both of $k,l$ differ from n (sf sequence $A145869$ in $[7]$):
$$2, 2, 4, 8, 30, 16, 28, 24, 114, 40, 58, 48, 100, 72, 92, 88,... .$$
On the other hand, for some even numbers such that
$$2, 4, 8, 16, 24, 28, 30, 32, 40, 48, 56, 58, 60, 62, 64, 72, 80, 88, 92, 96, 100, ...$$
(sf sequence $A145866$ in $[7]$) we have the uniqueness of representation (36).
\begin{question}\label{3}
What is the density  of such numbers?
\end{question}
\section{One-to-one mapping to lattice squares}
Notice that,
\begin{equation}\label{37}
s^{(r)}_{r^t+1}/s^{(r)}_{r^t}=r+1.
\end{equation}
Indeed, by Theorem 3,
\begin{equation}\label{38}
s^{(r)}_{r^t+1}=r^{2t+1}+1,
\end{equation}
while, since $r^{t}-1=(r-1)(1+r+...+r^{t-1}),$ then
\begin{equation}\label{39}
s^{(r)}_{r^t}=1+(r-1)(r+r^3+...+r^{2t-1})=\frac {r^{2t+1}+1} {r+1}.
\end{equation}
Therefore, the first number $N\equiv1 \pmod r$ which contains in its decomposition $s^{(r)}_{r^t+1}$ is $N=s^{(r)}_{r^t+1}+r,$ while every $N\leq s^{(r)}_{r^t+1}$ could contain in its decomposition only terms
$$s^{(r)}_{1}, s^{(r)}_{2},...,s^{(r)}_{r^t}.$$
 In view of (38), the number of such integers $N\geq r+1$ is
$$\lfloor\frac {r^{2t+1}+1} {r}\rfloor=(r^{t})^2.$$
Therefore, using the "coordinates" $k=k(r,N)$,  $l=l(r,N)$ of $N$, by Theorem 3, we have a one-to-one mapping
\begin{equation}\label{40}
\{N: N\equiv 1 \pmod r,\enskip r+1\leq N \leq r^{2t+1}+1\}\longleftrightarrow [1, r^{t}]\times[1, r^{t}],
\end{equation}\newpage
where $[1,n]=\{1,2,...,n\}.$
\begin{question}\label{4}
Let $\pi=\{\pi_1,...,\pi_{r^{2t}}\}$ be a permutation of numbers $r+1, 2r+1,...,r^{2t+1}+1.$ If to consider
the lattice square $[1, r^{t}]\times[1, r^{t}]$ as a planar graph in which the vertex corresponding to $\pi_i$ is
connected to the vertex corresponding to $\pi_{i+1},\enskip i=1,...,r^{2t}-1.$ What permutation $\pi$ corresponds
to minimal (maximal) Euclidean distance of the planar traveling salesman problem on the considered lattice square?
\end{question}
For example, in case of $r=2, t=1$ the permutation $\{3,5,9,7\}$  $(\{7,5,3,9\})$ corresponds to the minimal (maximal)
distance.

\;\;\;\;\;\;\;

\bfseries Acknowledgement.\enskip\mdseries The author is grateful to Professor Javier Cilleruelo (Madrid)
for his helpful and informative private correspondence respectively references [6] and [4].

\end{document}